\documentclass[reqno,11pt]{amsart}
\usepackage{amsmath, amssymb, amsthm,amsfonts} 
\usepackage[english]{babel}
\usepackage{bbm}
\usepackage{graphicx}
\usepackage{url}
\usepackage{epstopdf}
\usepackage[a4paper,bindingoffset=0.5cm,left=2cm,right=2cm,top=2.5cm,bottom=2cm,footskip=.8cm]{geometry}

\usepackage{tikz}
\usetikzlibrary{arrows, automata,positioning,calc,shapes,decorations.pathreplacing,decorations.markings,shapes.misc,petri,topaths}
\usepackage{tkz-berge}
  \usepackage{pgfplots}
  \pgfplotsset{compat=newest}
  \usetikzlibrary{plotmarks}
  \usepackage{grffile}
\newlength\figureheight
  \newlength\figurewidth
\pgfplotsset{%
    tick label style={font=\scriptsize},
    label style={font=\footnotesize},
    legend style={font=\footnotesize},
         every axis plot/.append style={very thick}
}

\usepackage{rotating}
\usepackage{amsbsy,enumerate}
\usepackage{graphicx}
\usepackage{ccaption}
\usepackage{comment}
\usepackage{mathrsfs}

\newcommand{\s}{^\star}

\newcommand{\bs}{\boldsymbol}

\newcommand{\vb}{\vspace{3.2mm}}

\newcommand{\vt}{\vartheta}
\newcommand{\DT}{\Delta t}

\allowdisplaybreaks

\newtheorem{lemma}{Lemma}

\newtheorem{remark}{Remark}

\newtheorem{proposition}{Proposition}

\makeatletter
\renewcommand{\fnum@figure}[1]{\textbf{\figurename~\thefigure}. }
\renewcommand{\fnum@table}[1]{\textbf{\tablename~\thetable}. }
\makeatother

\begin{document}

\title[A network of infinite-server queues with multiplicative transitions]{Networks of infinite-server queues\\ with multiplicative transitions}

\author{Dieter Fiems, Michel Mandjes, and Brendan Patch}

\begin{abstract}
This paper considers a network of infinite-server queues with the special feature that, triggered 
by specific events, the network population vector may undergo a linear transformation (a `multiplicative transition'). For this model we characterize the joint probability generating function in terms of a system of partial differential equations; this system
enables the evaluation of (transient as well as stationary) moments. We show that several relevant systems fit in the framework developed, such as networks of retrial queues, networks in which jobs can be rerouted when links fail, and storage systems. Numerical examples illustrate how our results can be used to support design problems. 

\vb

\noindent
{\sc Keywords.} Queueing networks $\circ$ infinite-server systems $\circ$ multiplicative transitions $\circ$ retrial queues $\circ$ rerouting $\circ$ storage systems

\vb

\noindent
{\sc Affiliations.} D.\ Fiems is with Department of Telecommunications and Information Processing, Ghent University, Ghent, Belgium. M.\ Mandjes and B.\ Patch are with Korteweg-de Vries Institute, University of Amsterdam; B.\ Patch is also with the School of Mathematics and Physics, The University of Queensland, St.\ Lucia, Australia. The research for this paper is partly funded by the NWO Gravitation Programme N{\sc etworks}, Grant Number 024.002.003 (Mandjes, Patch), an NWO Top Grant, Grant Number 613.001.352 (Mandjes), the ARC Centre of Excellence for Mathematical and Statistical Frontiers (Patch), and an Australian Government Research Training Program (RTP) scholarship (Patch).

\noindent Corresponding author: B.\ Patch. Address:  Korteweg-de Vries Institute for Mathematics, P.O. Box 94248,
1090 GE Amsterdam,
The Netherlands. Email: {\tt \footnotesize brendan.patch@gmail.com}.
\end{abstract}

\maketitle

\section{Introduction}

The vast majority of  queueing network models studied in the literature are of the following form: there is a set of nodes that are fed by  streams of external arrivals, and a routing mechanism that determines to which queue served clients are forwarded or whether the client leaves the system altogether.  The most common queueing disciplines are of single-server type (entailing that clients may have to wait until they get into service) and of infinite-server type (in which all customers present at a node are served in parallel).   

A key feature of  the conventional class of models  described above is that clients join and leave queues {\it one by one.}
In many applications, however, triggered by specific events, the {\it  full population} of individual queues may move around the network (or leave the system altogether). Particularly in the reliability and availability context, there are many relevant examples of such dynamics. We could for instance think of a data communication network with unreliable nodes: at the moment that a node goes down, all traffic residing at the node may be instantly lost.  Another example concerns rerouting: triggered, for instance, by a link failure, clients are moved from one set of resources to an alternative set (the `backup route'). Due to the fact that they correspond to transitions of the entire population of specific queues, the dynamics of the above two examples do not align with those of   conventional queueing models. 

\vb

{\it Scope, object of study.} Motivated by the above examples, the main objective of the present paper is to analyze queueing networks {\it with multiplicative transitions}. These multiplicative transitions effectively entail that the network dynamics include transitions by which the network's population vector, say ${\boldsymbol m}$, is (pre-)multiplied by a matrix $A$ with integer-valued, nonnegative entries, so that the network population after the transition becomes $A{\boldsymbol m}.$  For instance, choosing  $A$ to be a diagonal matrix with $[A]_{ii}=0$ and $[A]_{kk}=1$ for all $k\not= i$  would correspond to the event of all clients at node $i$ being lost. Relocation of clients can be taken care of in a similar manner: the full population of queue $i$ moving to queue $j$ corresponds to $[A]_{ji}=1$, $[A]_{kk}=1$ for all $k\not=i$, and all other entries equal to $0$. 

In this paper the queues considered are of {\it infinite-server} type. This type of queue  is particularly relevant in contexts where the sojourn time at a node  of each client is not (or hardly) affected by other  clients. As such, the model has a broad variety of applications, ranging from the number of websurfers simultaneously present at a set of websites, to the number of messenger RNA molecules simultaneously present in a collection of cells. A specific application that features in the present article concerns the optimal design of storage networks. To make the model as widely applicable as possible, we assume that all relevant transition rates (i.e., arrival rates and departure rates) are affected by an external autonomously evolving Markovian environmental process; the resulting model is therefore of a {\it Markov modulated} nature. As will become clear, in a reliability context such an environmental process can be used to model the state of the nodes of the network (i.e., `up' or `down'). 

\vb

{\it Contributions.} The paper has two main contributions. (i) In the first place we set up a general model of a network of infinite-server queues with multiplicative transitions. For this model we derive a system of partial differential equations that describe its time-dependent behavior (in terms of the probability generating function of the joint queue length distribution), as well as a procedure to evaluate the corresponding moments.  The model turns out to have a non-trivial stability condition (under which the system's stationary behavior is well-defined), which we establish using the expression we found for the time-dependent mean.  (ii) In the second place, we point out that various natural, practically relevant models fit in our framework. Most notably, we concentrate on a network of retrial queues, a network with rerouting, and a storage network. 
Our results can be used to support various design decisions. In 
the storage system application, for instance, interesting tradeoffs can be numerically assessed: files are typically stored on multiple locations so as to mitigate the risk of loss, but evidently one wants to do so without  using an unnecessarily large amount of storage space.

\vb

{\it Literature.} As mentioned above, in typical queueing network models the number of clients per queue changes by one at a time; see e.g.\ the standard textbooks \cite{KEL, SER}. Several papers, such as \cite{CHT,HT,MRS}, consider queues with batch arrivals and batch services and find product-form results, but these typically neither cover our multiplicative update rule nor allow the transition rates to be affected by an environmental process.

As mentioned above, a relevant special case of our model corresponds to the context of reliability. In many situations, when a network resource (a node or a link) fails, all clients using it will be lost. Such models are known as queueing models with {\it catastrophes}; for a fairly complete account of such models, we refer to the recent literature review in \cite[Section 1]{KPR}. The models studied are typically (but not always) one-dimensional;  interesting contributions include \cite{EF,SW}.

Queueing models for which the underlying infrastructure alternates between being `up' and `down' can be seen as  examples of stochastic processes on  dynamically evolving graphs. Despite the sizeable literature on random graphs, the body of work on dynamic random graphs is considerably smaller, and (evidently) the body of work covering stochastic processes on dynamic random graphs is even smaller. In a recent contributions, results on dynamic random graphs have been reported; see e.g. \cite{HO1,HO2,MSBS,ZMN}. Our paper can be regarded as being among the first to facilitate describing queueing processes on a randomly evolving graph (but it is noticed that the model we propose is substantially more general, as the multiplicative transitions are not restricted to node failures and repairs). 

As mentioned, our model covers various practically relevant models as special cases. In each of the corresponding application areas there is a large collection of papers and textbooks available; in Section 4 we include a number of domain-specific references.

\vb

{\it Organization.} The paper is organized as follows. Section \ref{migr} presents the model in its generic form, and after some preliminaries, the results in terms of partial differential equations characterizing the joint probability generating function and ordinary linear differential equations characterizing the moments. In addition, the stability condition is provided, under which stationary moments exist, which can be found by solving systems of linear equations. Section 3 gives an indication of the width of our framework: we show that it covers  a network of retrial queues, a network with rerouting, and a storage network. Section 4 demonstrates a couple of design issues that can be resolved by using our machinery. Finally, Section 5 provides a discussion and concluding remarks.

\section{Analysis}\label{migr}
This section studies our generic model: a network of infinite-server queues with  multiplicative transitions. We first introduce the model, then study its time-dependent behavior,  derive its stability condition,
and conclude by commenting on numerical issues. 
\subsection{Model}
In this subsection we describe our network of infinite-server queues with multiplicative transitions between the nodes. 
Let ${\mathscr N}:=\{1,\ldots,N\}$ (with $N\in{\mathbb N}$) be the set of infinite-server queues. The object of study is $({\bs M}(t))_{t\geqslant 0}$ (with ${\bs M}(t)\in{\mathbb N}_0^N$), that is, the multivariate queue content process (also sometimes referred to as the network population process). The process $(X(t))_{t\geqslant 0}$ (with $X(t)\in{\mathscr I}:=\{1,\ldots,I\}$) is the environment process (or: background process), which evolves autonomously of the queue content process; in our setup, $(X(t))_{t\geqslant 0}$ is assumed to be an irreducible continuous-time Markov chain.

The following transition rates play a role:
\begin{itemize}
\item[$\circ$]
The rate $\lambda_n^{(i)}\geqslant 0$ is the external arrival rate at queue $n\in\{1,\ldots,N\}$ when the background process $X(\cdot)$ is in $i \in {\mathscr I}$. Note that this entails that the arrival process at each of the queues is a Markov-modulated Poisson process. 
\item[$\circ$]
Likewise, the rate $\mu_{nn'}^{(i)}\geqslant 0$ is the departure rate {\it for every customer present} from queue $n$ to queue $n'$ when the background process $X(\cdot)$ is in $i \in {\mathscr I}$. Here $n\in{\mathscr N}$ and
$n'\in{\mathscr N}\cup\{0\}$, where $n'=0$ corresponds to the client leaving the network. Note that at the queues, the clients are served simultaneously, reflecting the infinite-server nature of each of the queues. 
\item[$\circ$]
Define for each pair $(i,j)$ with $i,j\in{\mathscr I}$ such that $i\not=j$ the set ${\mathscr K}_{ij}:=\{1,\ldots,K_{ij}\}$ with $K_{ij}\in{\mathbb N}.$
Let, for each $k\in {\mathscr K}_{ij}$,  $A_{ij}^{(k)}$ be an $(N\times N)$-matrix with entries in ${\mathbb N}_0$. 
The rate $\alpha_{ij}^{(k)}\geqslant 0$ is the rate at which the queue content, say ${\bs m}\in{\mathbb N}_0^N$,
is converted into $A_{ij}^{(k)} {\bs m}$, and at the same time the environment process $X(\cdot)$ jumps from state $i$ to state $j$, for $i,j\in {\mathscr I}$. For obvious reasons, we refer to these events as   {\it multiplicative transitions}. 
\end{itemize}
Two issues are worth highlighting. (i)~Note that the above description does not explicitly include notation for state transitions of the background process that do not involve multiplication with an $A$-matrix. It is easily observed, however, that such transitions can be introduced by letting the $A$-matrix correspond to an  identity matrix. (ii)~Transitions from $i$ to $j$ with  $i=j$ (`self-transitions') are allowed. Our setup in this respect differs from how continuous-time Markov processes are typically  defined; observe that $X(\cdot)$ is nonetheless a continuous-time Markov chain.

Notice that it can be anticipated that this system has a non-trivial stability condition. Observe that if some of the $A$-matrices have entries larger than 1, the parameters may be such that the network population can grow quickly and eventually explode. When the stability condition applies, however, this cannot happen. We derive the stability condition in Section \ref{stab}. 
Evidently, the system's time-dependent behavior can be studied regardless of the validity of 
such a stability condition; this time-dependent behavior is the topic of Sections \ref{spde}--\ref{mom}. In Section \ref{NUMER} we comment on the numerical evaluation of the performance measures under study. 

\subsection{System of partial differential equations}\label{spde}
The objective of this subsection is  to characterize the distribution of $({\bs M}(t),X(t))\in{\mathbb N}_0^N\times {\mathscr I}$. We take the classical approach of setting up a system of partial differential equations for the corresponding transforms. To this end, we first define, for $i\in{\mathscr I}$ and $t\geqslant 0$, 
\[\varphi_i({\bs \vt},t) :={\mathbb E}\big({\rm e}^{\langle {\bs \vt},{\bs M}(t)\rangle}I_i(t)\big),\]
with $I_i(t) :=1\{X(t)=i\}$, the indictor function for the event that $X(t)$ equals $i$. Evidently, these quantities uniquely describe the system's probabilistic behavior.

So as to set up the differential equations, the main idea is to relate the state of the system at time $t+\DT$ to the state at time $t$, for $\DT$ small.
We rely  on the usual `Markovian reasoning', meaning that when the environmental process is in state $i$ at time $t$ the following three types of events have to be considered:  (i)~with a probability of essentially $\lambda_n^{(i)}\DT$ there is an external arrival at node $n$, 
(ii)~with probability $\mu_{nn'}^{(i)}{M}_n(t)\,\DT$ there is a departure from node $n$ to $n'$ (with $n'$ possibly equalling $0$, to model the clients that leave the network), and (iii)~with probability $\alpha_{ij}^{(k)}\DT$ the environmental process jumps to $j$ while simultaneously the network population vector ${\bs M}(t)$ is instantly replaced by $A_{ij}^{(k)}{\bs M}(t)$. Working out these transitions in detail, 
elementary calculations reveal that, as $\DT\downarrow 0$,
\begin{equation}\label{de}
\varphi_i({\bs \vt},t+\DT) =\varphi_i({\bs \vt},t)+\xi_i({\bs\vt},t) \DT -\zeta_i({\bs \vt},t)\DT+{\rm o}(\DT),\end{equation}
where
\begin{eqnarray*}
\xi_i({\bs \vt},t) &:=&\sum_{n=1}^N {\rm e}^{\vt_n} {\mathbb E}\big({\rm e}^{\langle {\bs \vt},{\bs M}(t)\rangle}I_i(t)\big) \,\lambda_n^{(i)} +
 \sum_{n=1}^N\sum_{n'=1}^N {\rm e}^{-\vt_n+\vt_{n'}}
{\mathbb E}\big({\rm e}^{\langle {\bs \vt},{\bs M}(t)\rangle}I_i(t)\,M_n(t)\big)\,\mu_{nn'}^{(i)} +\\
&&\sum_{n=1}^N {\rm e}^{-\vt_n}
{\mathbb E}\big({\rm e}^{\langle {\bs \vt},{\bs M}(t)\rangle}I_i(t)\,M_n(t)\big) \,\mu_{n0}^{(i)} +
\sum_{j=1}^I  \sum_{k=1}^{K_{ji}} {\mathbb E}\big({\rm e}^{\langle {\bs \vt},A_{ji}^{(k)}{\bs M}(t)\rangle}I_j(t)\big) \,\alpha_{ji}^{(k)}
\end{eqnarray*}
and
\begin{eqnarray*}
\zeta_i({\bs \vt},t) &:=&\sum_{n=1}^N  {\mathbb E}\big({\rm e}^{\langle {\bs \vt},{\bs M}(t)\rangle}I_i(t)\big)  \,\lambda_n^{(i)} + \sum_{n=1}^N\sum_{n'=1}^N {\mathbb E}\big({\rm e}^{\langle {\bs \vt},{\bs M}(t)\rangle}I_i(t)\,M_n(t)\big) \,\mu_{nn'}^{(i)} +\\
&&\sum_{n=1}^N {\mathbb E}\big({\rm e}^{\langle {\bs \vt},{\bs M}(t)\rangle}I_i(t)\,M_n(t)\big) \,\mu_{n0}^{(i)} +
\sum_{j=1}^I\sum_{k=1}^{K_{ij}} {\mathbb E}\big({\rm e}^{\langle {\bs \vt},{\bs M}(t)\rangle}I_i(t)\big) \alpha_{ij}^{(k)} .
\end{eqnarray*}
To understand the structure of $\xi_i({\bs \vt},t)$ and $\zeta_i({\bs \vt},t)$, note that 
their first terms reflect the external arrivals, the second terms the routing to other queues, the third terms clients leaving the network, and the fourth terms the multiplicative transitions. 

The next step is to rewrite the expressions for $\xi_i({\bs \vt},t)$ and $\zeta_i({\bs \vt},t)$ in terms of partial derivatives with respect to the arguments $\vartheta_n$, $n\in{\mathscr N}$. We thus obtain, with $A^{\rm T}$ denoting the transpose of the matrix $A$,
\begin{eqnarray*}
\xi_i({\bs \vt},t) &=&
\sum_{n=1}^N  {\rm e}^{\vt_n} \varphi_i({\bs\vt},t)\,\lambda_n^{(i)} +
 \sum_{n=1}^N\sum_{n'=1}^N {\rm e}^{-\vt_n+\vt_{n'}}\frac{\partial}{\partial \vt_n}\varphi_i( {\bs \vt},t) \,\mu_{nn'}^{(i)} +\\
&&\sum_{n=1}^N {\rm e}^{-\vt_n}\frac{\partial}{\partial \vt_n}\varphi_i( {\bs \vt},t)  \,\mu_{n0}^{(i)} +\sum_{j=1}^I \sum_{k=1}^{K_{ji}} 
 \varphi_j\big((A_{ji}^{(k)})^{\rm T}{\bs\vt},t\big) \alpha_{ji}^{(k)} 
 \end{eqnarray*}
and
\begin{eqnarray*}
\zeta_i({\bs \vt},t) &=&
\sum_{n=1}^N   \varphi_i({\bs\vt},t)\,\lambda_n^{(i)} +
 \sum_{n=1}^N\sum_{n'=1}^N \frac{\partial}{\partial \vt_n}\varphi_i( {\bs \vt},t) \,\mu_{nn'}^{(i)} +\\
&&\sum_{n=1}^N \frac{\partial}{\partial \vt_n}\varphi_i( {\bs \vt},t)  \,\mu_{n0}^{(i)} +\sum_{j=1}^I \sum_{k=1}^{K_{ij}} 
 \varphi_i({\bs\vt},t)
 \alpha_{ij}^{(k)} .
 \end{eqnarray*}
 We proceed in the common way: by subtracting $\varphi_i( {\bs \vt},t)$ from both sides in (\ref{de}), dividing by $\DT$, and taking the limit $\DT\downarrow 0$, we arrive at the following result. \begin{proposition}
 The transforms $\varphi_i({\bs \vt},t)$, for $i\in{\mathscr I}$, satisfy the following system of partial differential equations:
 \begin{eqnarray}
 \frac{\partial}{\partial t}\varphi_i( {\bs \vt},t)&=&
\varphi_i({\bs\vt},t)\sum_{n=1}^N \big( {\rm e}^{\vt_n} -1\big)\,\lambda_n^{(i)} +
 \sum_{n=1}^N\sum_{n'=1}^N\big( {\rm e}^{-\vt_n+\vt_{n'}}-1\big)\frac{\partial}{\partial \vt_n}\varphi_i( {\bs \vt},t) \,\mu_{nn'}^{(i)} +\notag\\&&
 \sum_{n=1}^N \big({\rm e}^{-\vt_n}-1\big)\frac{\partial}{\partial \vt_n}\varphi_i( {\bs \vt},t)  \,\mu_{n0}^{(i)} +\notag\\&&
 \sum_{j=1}^I \sum_{k=1}^{K_{ji}}   \varphi_j\big((A_{ji}^{(k)})^{\rm T}{\bs\vt},t\big) \alpha_{ji}^{(k)}
-  \varphi_i( {\bs \vt},t) \sum_{j=1}^I \sum_{k=1}^{K_{ij}}  \alpha_{ij}^{(k)} \, . \label{mgf}
  \end{eqnarray}
  \end{proposition}
  From this relation moments can be evaluated by differentiation and inserting ${\boldsymbol\vartheta}={\boldsymbol 0}$, as we demonstrate in the next subsection.

\subsection{Moments}\label{mom}  
We now find an explicit expression for the $I$ mean queue content vectors (each of them in ${\mathbb R}_+^N$)
\[\bar {\bs M}_i(t)=\big[{\mathbb E}(M_n(t) I_i(t))\big]_{n=1}^N,\] $i \in \mathscr I$. 
In addition to playing a central role in our performance evaluation framework, these expressions also allow us  to establish the stability condition for this type of queueing network; see Section \ref{stab}.

The first step is to identify the transient distribution of the environmental process $X(\cdot)$. To this end, we let $\pi_i(t):= \mathbb P(X(t)=i)$ for $i \in \mathscr I$; this means that $\bs\pi(t) = [\pi_i(t)]_{i=1}^I$ is the transient distribution vector of the background process. Inserting $\bs \vt = \bs 0$ in (\ref{mgf})  yields a (homogeneous)  system of coupled linear differential equations:
 \begin{eqnarray*}
 \pi_i'(t)&=& \sum_{j=1}^I \sum_{k=1}^{K_{ji}}  \pi_j(t)  \alpha_{ji}^{(k)}-\pi_i(t) \sum_{j=1}^I \sum_{k=1}^{K_{ij}}    \alpha_{ij}^{(k)} .
  \end{eqnarray*}
This system can be compactly rewritten as  
\[
\bs \pi'(t) = \bar{\mathscr A}^{\rm T} \bs\pi(t)
\]
with $\bar{\mathscr A}=[\bar\alpha_{ij}]_{i,j=1}^I$ the corresponding generator matrix with elements,
for $i,j \in \mathscr I$ and $i\ne j$,
\[
\bar\alpha_{ij} = \sum_{k=1}^{K_{ij}} \alpha_{ij}^{(k)} \, , \quad \bar\alpha_{ii} = -\sum_{i'\not=i} \bar\alpha_{ii'} \, .
\]
We thus find ${\bs\pi}(t) =\exp(\bar{\mathscr A}^{\rm \,T} t){\bs\pi}(0).$ Observe that $\bar{\mathscr A}$ is a transition rate matrix (i.e., a matrix with  negative diagonal elements  and row sums equal to zero). This entails that, for any $t\geqslant 0$, ${\bs\pi}(t)$ is a probability distribution on ${\mathscr I}$.

Our next objective is to identify the first moments of the queue sizes. 
To obtain these quantities we differentiate the full equation (\ref{mgf}) with respect to each of the $\vartheta_n$ ($n\in{\mathscr N}$). Plugging in $\bs \vt = \bs 0$ then leads, after some straightforward but tedious algebra, to the following (non-homogeneous) system of 
linear differential equations:
 \[\bar {\bs M}_i'(t)
=
{\mathscr L}_i    \pi_i(t)+ {\mathscr M}_i \bar {\bs M}_i(t)
+ \sum_{j=1}^I {\mathscr A}_{ji} \bar {\bs M}_j(t) ,
  \]
  with the matrices ${\mathscr L}_i$, ${\mathscr M}_i$, and ${\mathscr A}_{ij}$ defined as follows.
    \begin{itemize}
  \item[$\circ$]
Firstly,  ${\mathscr L}_i := \big[\lambda_n^{(i)}\big]_{n=1}^N$, i.e., a column vector  with   the arrival rates in the different queues when the background process is in state $i\in{\mathscr I}$. 
\item[$\circ$]Secondly,
 \[{\mathscr M}_i :=\big[\mu_{n'n}^{(i)}+1\{n=n'\}\bar\mu_n^{(i)}\big]_{n,n'=1}^N,\:\:\:\mbox{with}\:\:\bar\mu_{n}^{(i)} = -\sum_{n'=0}^N \mu_{nn'}^{(i)},\] is the matrix with the departure rates between the different queues when the background process is in state $i\in{\mathscr I}$. 
 \item[$\circ$]In addition, we introduce the following notation for the multiplicative update process:
\[
{\mathscr A}_{ij} := \sum_{k=1}^{K_{ij}}  \alpha_{ij}^{(k)} \,A_{ij}^{(k)}  \, , \quad {\mathscr A}_{ii} := \sum_{k=1}^{K_{ii}}  \alpha_{ii}^{(k)} \,A_{ii}^{(k)}  -  \bar\alpha_{i} {\mathbb I}_N \, ,
\]
for $i,j \in \mathscr I$, $i \ne j$, with ${\mathbb I}_N$ denoting the $(N\times N$)-dimensional identity matrix, and with
\[\bar\alpha_i:= \sum_{j=1}^I \sum_{k=1}^{K_{ij}}\alpha_{ij}^{(k)}.\]
\end{itemize}

Moreover, the above  set of equations can be a combined into a single equation. To this end, we define  $\bar {\bs M}(t)$ to be the vector $ [\bar {\bs M_i}(t)]_{i=1}^I$ of dimension $J:=IN$. Also
${\mathscr A} := [{\mathscr A}_{ji}]_{i,j=1}^I$ and ${\mathscr M} := \text{diag}([{\mathscr M}_i]_{i=1}^I)$, which are  $(J\times J)$-dimensional matrices. Finally, ${\mathscr L} := \text{diag}([{\mathscr L}_i]_{i=1}^I)$ is a matrix of dimension $J\times I$. 
\begin{proposition}\label{prop:mean}
For any $t\geqslant 0$,
 \begin{align}
\bar {\bs M}'(t)
&=
 {\mathscr L} {\bs \pi}(t)  + ({\mathscr M} + {\mathscr A}) \bar {\bs M}(t) \, . \label{mom1}
  \end{align}
\end{proposition}

Solving the differential equation for the transient moment vector \eqref{mom1} leads to the following explicit solution (in terms of integrals over matrix-exponentials):
\begin{align}\nonumber
\bar {\bs M}(t) &= {\rm e}^{({\mathscr M}+{\mathscr A})t} \bar {\bs M}(0)  + \int_0^t   {\rm e}^{({\mathscr M}+{\mathscr A})(t-s)} {\mathscr L} {\bs \pi}(s) \, {\rm d}s \\
&=
{\rm e}^{({\mathscr M}+{\mathscr A})t} \bar {\bs M}(0)  + \int_0^t   {\rm e}^{({\mathscr M}+{\mathscr A})(t-s)} {\mathscr L} \,{\rm e}^{\bar{\mathscr A}^{\rm \,T} s}{\bs\pi}(0)
 \, {\rm d}s.\label{MEM}
\end{align}
The stationary means follow from equating $\bar {\bs M}(t)$ to ${\bs 0}$  and defining 
${\bs \pi}$ as the solution to ${{\mathscr A}}^{\rm T}{\bs\pi} = {\bs 0}$,  so that the stationary mean $\bar {\bs M}$ is given by
\begin{equation}\label{statm}
 \bar {\bs M} = -({\mathscr M}+{\mathscr A})^{-1} {\mathscr L}  {\bs \pi} 
.\end{equation}
Note, however, that this reasoning tacitly assumes that the underlying queueing network is stable, an issue we return to   in Section \ref{stab}.

Along the same lines higher  moments of the queue sizes can   be found as well. The higher transient moments can be phrased in terms of a (non-homogeneous) system of
linear differential equations. The procedure to identify them is of a recursive nature, as determining the $n$-th transient moment requires knowledge of the first   $n{-}1$
transient moments. Similarly, higher stationary moments follow as solutions to linear equations (under the stability condition), where finding the $n$-th stationary moment
requires the first $n{-}1$ 
stationary moments being available.  For analogous procedures in a related context, see \cite{KOEN}.

\subsection{Stability}\label{stab}
As it turns out, Prop.\ \ref{prop:mean} facilitates the  provision of conditions for the ergodicity of the Markov chain. Before proceeding to stating and proving our stability result, we first define $\omega$ to be spectral abscissa  of ${\mathscr M}+{\mathscr A}$, that is
\[
\omega := \max\{{\rm R}{\rm  e}\,\lambda: \lambda \in \text{spec}({\mathscr M}+{\mathscr A}) \}
\]
where $\text{spec}({\mathscr M}+{\mathscr A})$ is the set of eigenvalues of ${\mathscr M}+{\mathscr A}$. 

\begin{proposition}
The Markov chain $\bs Z(t) \equiv (\bs M(t),X(t)) $ is ergodic provided  $\omega$  is negative.
\end{proposition}

{\it Proof.} To prove the claim,   we study the ergodicity of the {\it skeleton Markov chain} $\{\bs Z(\Delta n); n \in \mathbb N\}$ for some $\Delta > 0$. Obviously, if the skeleton Markov chain is ergodic for some $\Delta > 0$, so is $\bs Z(t)$, as the mean recurrence time for any state of the skeleton chain is an upper bound for the mean recurrence time of the original chain $\bs Z(t)$. 

Appealing to \cite[Prop.\ I.5.3]{AsmussenAPQ}, it suffices to show that for  some $\varepsilon > 0$, $\Delta >0$, and $\mathfrak M \geqslant 0$,  the following mean drift condition holds:
\[
\mathbb E_{(\bs m,i)}\big(\Vert \bs M(\Delta)\Vert_1\big)-\Vert \bs m\Vert_1 < -\varepsilon
\]
for all $\bs m$ with $\Vert \bs m\Vert_1 > \mathfrak M$ and all $i \in \mathscr I$; the subscript $(\bs m,i)$ indicates that  the expectation is conditional on $ \bs Z(0) = (\bs m, i)$.

Define
$L:= \sum_{i=1}^I\sum_{n=1}^N \lambda_n^{(i)}$.
From the differential equation for the first moment \eqref{mom1}, we derive the following bound:\begin{align*}
\mathbb E_{(\bs m,i)}\big(\Vert \bs M(\Delta)\Vert_1\big)-\Vert \bs m\Vert_1 
&= \left\Vert \bar {\bs M}(\Delta)|_{\bs M(0)=\bs m, \bs \pi(0) = \bs e_i}\right\Vert_1  - \Vert \bs m\Vert_1\\
&= \left\Vert {\rm e}^{({\mathscr M}+{\mathscr A})\Delta} (\bs e_i \otimes \bs m)  + \int_0^\Delta  {\rm e}^{({\mathscr M}+{\mathscr A})(\Delta-s)} {\mathscr L} {\bs \pi}(s) {\rm d}s \right\Vert_1 - \Vert \bs m\Vert_1 \\
&\leqslant \Vert {\rm e}^{({\mathscr M}+{\mathscr A})\Delta}  \Vert_1 \Vert \bs m  \Vert_1 + \int_0^\Delta  \Vert {\rm e}^{({\mathscr M}+{\mathscr A})(\Delta-s)} \Vert_1\,L\, {\rm d}s - \Vert \bs m\Vert_1  \\
&\leqslant \Vert\bs m\Vert_1 \gamma {\rm e}^{\omega^{\star} \Delta} + \gamma  \,L \int_0^\Delta   {\rm e}^{\omega^{\star} (\Delta-s)} {\rm d}s  - \Vert \bs m\Vert_1\\
&\leqslant  \Vert\bs m\Vert_1  \gamma {\rm e}^{\omega^{\star} \Delta} - \frac{\gamma}{\omega^{\star}} \,L\,  (1 - {\rm e}^{\omega^{\star} \Delta}) -  \Vert\bs m\Vert_1  =: g(\Delta) \, ,
\end{align*}
for $\omega <  \omega^{\star} < 0$ and where $\gamma > 0$ is a constant; see \cite[Prop.\  11.18.8]{Bernstein} for the bound on the norm of the matrix exponential. Clearly, with $R:= \gamma L  /\omega^{\star}<0,$ 
\[\lim_{\Delta \to \infty} g(\Delta) = -R  - \Vert \bs m \Vert_1,\] which is negative for $\Vert \bs m \Vert_1 >  - R>0$. Hence, for any choice of $\mathfrak M > - R$, there exists a $\Delta$ such that the drift condition  of the skeleton chain holds for $\Vert \bs m \Vert_1 > \mathfrak M$. The skeleton chain $\{\bs Z(\Delta n)\}_{n \in \mathbb N}$ is therefore ergodic, which is inherited by  the original Markov process.\hfill$\Box$

\vb

\begin{remark} {\em At first sight it may  look unnatural that the stability condition is in terms of the matrices ${\mathscr M}$ and ${\mathscr A}$ only, and does not involve the external arrival rate matrix ${\mathscr L}$. To get a feel for the underlying intuition, let us consider the simplest  network possible: an isolated infinite-server queue, with external arrival rate $\lambda\geqslant 0$, exponential holding times with mean $\mu^{-1}\geqslant 0$, and a multiplicative transition from state $m\in{\mathbb N}_0$ to $a_k m$ (with $a_k\in{\mathbb N}_0$) with rate $\alpha_k\geqslant 0$ ($k=1,\ldots,K$). Then, using the  results of Section \ref{mom}, the mean number of clients in the queue at time $t$, denoted by $\bar M(t)$, satisfies the differential equation
\[\bar M'(t) = \lambda + \sum_{k=1}^K \alpha_k (a_k-1) \bar M(t) - \mu \bar M(t);\]
observe that the process goes up by one with rate $\lambda$, jumps from $m$ to $a_k m$ (leading to a net change of $(a_k-1)m$) with rate $\alpha_k$, and goes down by one with rate $\mu m$. 
To ensure stability, the mean number in the system should not explode.
This leads  to a stability condition that does not involve $\lambda$, viz.\ (in self-evident vector notation)
$\langle {\bs \alpha},{\bs a}-{\bs 1}\rangle <\mu.$} \hfill$\Diamond$\end{remark}

\subsection{Efficient evaluation of performance metrics} \label{NUMER}
In many applications, the performance of the system during a finite time interval $[0,T]$ is expressed in terms of quantities of the form 
\[v(T):=\sum_{i=1}^I\sum_{n=1}^N \varrho_{n,i} \bar M_{n,i}(T),\:\:\:\:w(T):=\int_0^T\sum_{i=1}^I\sum_{n=1}^N  \varrho_{n,i} 
\bar M_{n,i}(t){\rm d}t,\]
for some vector ${\bs \varrho}\in {\mathbb R}^J$ and $T\geqslant 0$,  
with $\bar M_{n,i}(t):={\mathbb E} \big(M_n(t)I_i(t)\big)$. In this section we point out how to efficiently compute the vectors $\bar {\bs M}(T)$ and $\int_0^T \bar {\bs M}(t){\rm d}t.$ 

We first study  $\bar {\bs M}(T)$; note that  $v(T)$ then follows upon evaluation of $\langle{\bs\varrho},
\bar {\bs M}(T)\rangle$. The first term of  expression (\ref{MEM}) 
is a matrix-exponential, for which standard evaluation techniques have been developed; see e.g.\
\cite{MV}. The second term reads $B(T)\cdot {\bs \pi}(0)$, with
\[B(T):=\int_0^T {\rm e}^{({\mathscr M}+{\mathscr A}) (T-s)} {\mathscr L} {\rm e}^{\bar {\mathscr A}^{\rm\, T} s} {\rm d}s .\]
By   \cite[Thm.\ 1]{VL}, $B(T)$ equals the $(J\times I)$-dimensional top right corner matrix of ${\rm e}^{{\mathscr C}T}$, where
\[{\mathscr C}:=\left[\begin{array}{cc}{\mathscr M}+{\mathscr A}&{\mathscr L}\\{\mathbb O}_{I,J}&\bar {\mathscr A}^{\rm\, T}\end{array}\right]\]
(with ${\mathbb O}_{I,J}$ defined as an all-zeros matrix of dimension $I\times J$).  We thus end up with the
following result.
\begin{lemma} \label{lemma1}
For any $T\geqslant 0$,
\[\bar {\bs M}(T) ={\rm e}^{({\mathscr M}+{\mathscr A})T} \bar {\bs M}(0) +
\big[ {\mathbb I}_J,{\mathbb O}_{J,I}\big]\cdot
{\rm e}^{{\mathscr C}T}
\cdot
\left[
\begin{array}{c}{\mathbb O}_{J,I}\\{\mathbb I}_I
\end{array}\right]
{\bs\pi}(0)
.\]
 \end{lemma}
Now we explain how to evaluate  $\int_0^T\bar {\bs M}(t){\rm d}t$, which facilitates the computation of 
\[w(T) =\int_0^T \langle{\bs \varrho}, \bar {\bs M}(t)\rangle{\rm d}t = \left\langle{\bs\varrho},
\int_0^T\bar {\bs M}(t){\rm d}t\right\rangle.\] 
Due to Lemma \ref{lemma1},
\[ \int_0^T\bar {\bs M}(t){\rm d}t =
\int_0^\infty {\rm e}^{({\mathscr M}+{\mathscr A})t} {\rm d}t\cdot \bar {\bs M}(0) +
\big[ {\mathbb I}_J,{\mathbb O}_{J,I}\big]\cdot \int_0^T
{\rm e}^{{\mathscr C}t}{\rm d}t
\cdot
\left[
\begin{array}{c}{\mathbb O}_{J,J}\\{\mathbb I}_I
\end{array}\right]
{\bs\pi}(0)
.\]
Define, with $J^+:= J+I$, the matrices
\[ 
{\mathscr C}_1:=\left[\begin{array}{cc}
{\mathbb O}_{J,J}
&{\mathbb I}_{J}\\{\mathbb O}_{J,J}&{\mathscr M}+{\mathscr A}\end{array}\right],\:\:\:
{\mathscr C}_2:=\left[\begin{array}{cc}
{\mathbb O}_{J^+,J^+}
&{\mathbb I}_{J^+}\\{\mathbb O}_{J^+,J^+}&
{\mathscr C}\end{array}\right],\]
which are of dimensions $2J\times 2J$ and $2J^+\times 2J^+$, respectively.
Again applying \cite[Thm.\ 1]{VL}, we arrive at 
\begin{align*}
\int_0^T\bar {\bs M}(t){\rm d}t=&\big[ {\mathbb I}_J,{\mathbb O}_{J,J}\big]\cdot
{\rm e}^{{\mathscr C}_1T}
\cdot
\left[
\begin{array}{c}{\mathbb O}_{J,J}\\{\mathbb I}_J
\end{array}\right]
{\bs M}(0)\,+\\
&\big[ {\mathbb I}_J,{\mathbb O}_{J,I}\big]\cdot 
\big[ {\mathbb I}_{J^+},{\mathbb O}_{J^+,J^+}\big]\cdot
{\rm e}^{{\mathscr C}_2 T} 
\cdot
\left[
\begin{array}{c}{\mathbb O}_{J^+,J^+}\\{\mathbb I}_{J^+}
\end{array}\right]
\left[
\begin{array}{c}{\mathbb O}_{J,J}\\{\mathbb I}_I
\end{array}\right]
{\bs\pi}(0).
\end{align*}
This can be rewritten in the following more compact form. 
\begin{lemma} \label{lemma2}
For any $T\geqslant 0$,
\begin{align*}
\int_0^T\bar {\bs M}(t){\rm d}t=&\big[ {\mathbb I}_J,{\mathbb O}_{J,J}\big]\cdot
{\rm e}^{{\mathscr C}_1T}
\cdot
\left[
\begin{array}{c}{\mathbb O}_{J,J}\\{\mathbb I}_J
\end{array}\right]
{\bs M}(0)+\big[ {\mathbb I}_J,{\mathbb O}_{J,2J^+-J}\big]\cdot
{\rm e}^{{\mathscr C}_2 T} 
\cdot
\left[
\begin{array}{c}{\mathbb O}_{2J^+-I,I}\\{\mathbb I}_{I}
\end{array}\right]
{\bs\pi}(0).
\end{align*}
 \end{lemma}

\section{Retrial queues, rerouting, storage systems}
In this section we show the power of the framework introduced in the previous section, by 
pointing out how it facilitates the modelling of all sorts of relevant phenomena. We specifically focus on: (i)~systems in which nodes go down but in which lost customers attempt reentry, (ii)~systems in which customers are rerouted when one of the links along the route goes down, and (iii) storage systems. 

\subsection{Retrial queues}
In this subsection we consider a network of faulty service stations. Each of the stations alternates between being `up' and `down'. While a station is in the `up' state it processes clients as a standard infinite-server queue. Upon going down, all clients present at a service station move instantly to an associated retrial location, from where they (independently of each other) try to reenter the service station or renege. For an in-depth account of related retrial models, we refer to~\cite{JA}. We note  that, to the best of our knowledge, the literature does not cover the model we study here.

We now point out how this retrial mechanism fits in the framework that we set up in the previous section.
Let the components $1$ up to $N^\circ$ of ${\bs M}(\cdot)$ correspond to the service stations in the network, and the components $N^\circ+1$ up to $2N^\circ =: N$ to the associated retrial locations. 
Here we assume that the up-time of station $n\in\{1,\ldots,N^\circ\}$ is exponentially distributed with parameter $\gamma_n^{({\rm u})}$, and the corresponding down-time is exponentially distributed with parameter $\gamma_n^{({\rm d})}$. We thus have constructed an environmental process of dimension $I=2^{N^\circ}$, where each state of this process corresponds to the particular set of stations that are up (and consequently also the set of stations that are down). In the sequel we let $S(i)$ denote the set of stations that are up when the environmental process is in state $i$. (It is noted that the above model can be extended in a straightforward manner to the situation in which the up-times and down-times stem from phase-type distributions. Similarly, Markov-modulated arrivals can be dealt with.)

We let $\lambda_n$ be the arrival rate at station $n$;
note that clients arriving at station $n$ when it is down are immediately placed in the corresponding retrial pool (which is component $N^\circ +n$ of ${\bs M}(\cdot)$).
 Also, let $\mu_{nn'}$ denote the rate of being routed (after service completion) from node $n$ to node $n'$ (with $n'=0$ corresponding, as always, with the client leaving the network). The rate $\kappa_n$ is the retrial rate at the $n$-th retrial location (i.e., component $N^\circ +n$ of ${\bs M}(\cdot)$), and $\nu_n$ the corresponding renege rate (reflecting clients that leave the network from a retrial location, i.e., before being served, e.g.\ due to impatience). 

Let us now describe how the above parameters translate into the rates of the framework of the previous section. Suppose  the environmental process is in state $i$. Let us first consider the external arrivals. Define $1_n(i):=1\{n \in S(i)\}$.
For $n=1,\ldots,N^\circ$, the external arrival rates when the environmental process is in state  $i$, are given by
\[\lambda_n^{(i)} = \lambda_n\,1_n(i),\:\:\: \lambda^{(i)}_{n+N^\circ} = \lambda_n\,(1-1_n(i)).\]
Regarding the service completions, we have for the service stations (with $n,n'=1,\ldots,N^\circ$)
\[\mu_{nn'}^{(i)} = \mu_{nn'}\,1_{n'}(i),\:\:\:\:\mu_{n,n'+N^\circ}^{(i)}=
\mu_{nn'}\,(1-1_{n'}(i)),
\:\:\:\:\mu_{n0}^{(i)} = \mu_{n0},\]
and for the retrial locations (again with $n=1,\ldots,N^\circ$)
\[\mu^{(i)}_{n+N^\circ,n} = \kappa_n\,1_n(i),\:\:\:\mu^{(i)}_{n+N^\circ,0} =\nu_n.\]
We now consider the  transitions related to the stations alternating between the active and inactive mode. Two cases are to be distinguished; as it turns out, for all $i,j\in{\mathscr I}$ we have that $K_{ij}=1$. 
\begin{itemize}
\item[$\circ$] Suppose that for a $j\not= i$ and some $n\in\{1,\ldots,N^\circ\}$   we have $S(i)=S(j)\cup\{n\}$; then the background process jumps from $i$ to $j$ after an exponentially distributed time with rate $\alpha_{ij}=\gamma_n^{({\rm u})}$. Note that this transition corresponds to station $n$ failing, and thus clients being moved to the corresponding retrial location. The  vector ${\bs M}(t)$  is premultiplied by a $(N\times N)$-dimensional matrix $A_{ij}$ consisting of a $0$ on position $(n,n)$, a~$1$ on position $(n+N^\circ,n)$, all diagonal entries except the $n$-th being 1, and all other entries being~0.
\item[$\circ$] 
Suppose on the other hand that for $i\not=j$ and some $n\in\{1,\ldots,N^\circ\}$ we have $S(j)=S(i)\cup\{n\}$; then the background process jumps from $i$ to $j$ with rate $\alpha_{ij}=\gamma_n^{({\rm d})}$, without the network population vector changing. This transition corresponds to station $n$ having been repaired. 
\end{itemize}

\vb

This framework has the potential to support various design issues. In the network described, an objective may be to keep the fraction of lost clients (due to reneging) below some tolerable level, say, $\varepsilon$. 
To this end, define $Z_a(t)$ as the total number of clients arrived in $[0,t]$ and $Z_\ell(t)$ as the number of clients lost. With ${\bs \lambda}$ defined in the evident way,
\[{\mathbb E}\,Z_a(t) = \int_0^t 
\sum_{i=1}^I \sum_{n=1}^N \pi_i(s)\lambda_n^{(i)}{\rm d}s =
\int_0^t \langle  {\bs \lambda},\bar {\bs \pi}(s)\rangle {\rm d}s.\]
Likewise, with ${\bs \eta}$ defined appropriately (i.e., a vector of which the first $IN^\circ$ entries equals 0 and the  second $IN^\circ$ entries equal the appropriate $\nu_n$),
\begin{equation}
\label{LOST}{\mathbb E}\,Z_\ell(t) = \int_0^t 
\sum_{i=1}^I \sum_{n=N^\circ+1}^N {\mathbb E}\big(M_n(s)I_i(s) \big)\nu_n\, {\rm d}s =
\int_0^t \langle  {\bs \nu},\bar {\bs M}(s)\rangle {\rm d}s.\end{equation}
The numerical evaluation of the above performance metrics is facilitated by Lemma \ref{lemma2}.

Suppose that for a given time horizon $T$ the service requirement is ${\mathbb E}\,Z_\ell(T)\leqslant 
\varepsilon\cdot{\mathbb E}\,Z_a(T)
$. If for given repair rates ${\bs \gamma}^{({\rm d})}\equiv (\gamma_1^{({\rm d})},\ldots ,\gamma_{N^\circ}^{({\rm d})})$ this condition is not met, one may wonder by how much the repairs have to be sped up to meet the service requirement. A relevant optimization problem is then
\[\min_{{\bs \gamma}^{({\rm d})}}\, \langle {\bs \gamma}^{({\rm d})},{\bs 1}\rangle,\:\:\:\mbox{subject to}\:\:{\mathbb E}\,Z_\ell(T)\leqslant 
\varepsilon\cdot{\mathbb E}\,Z_a(T)
.\]

\subsection{Rerouting}
Routing concerns the selection of a path along which traffic is transmitted. To make the service level more robust,
one may adopt the policy that when a network element fails, traffic using that network element is routed along an alternative route.
For a textbook treatment of routing in communication networks, we refer to e.g.\ \cite{KR}.

Our present framework can be used to track the number of clients that use the different direct and indirect routes. The clients along these routes correspond to the customers of our framework and the queues (i.e., the components of ${\bs M}(\cdot)$) record the quantity of clients utilizing each of the  direct and indirect routes. More formally, the rerouting model can be cast in our framework as follows. 
Let there be $N^\circ$ origin-destination pairs, each connected by a direct route (consisting of one link)  as well as an indirect route (consisting of two links). Let the direct link used by the $n$-th origin-destination pair be labelled by $n$, and let ${\mathscr N}(n):=\{n_1(n),n_2(n)\}$ (both elements being contained in $\{1,\ldots,N^\circ\}\setminus \{n\}$) be the links of the corresponding indirect route. We thus have $N=2N^\circ$ queues, the first $N^\circ$ queues corresponding to the number of clients on the direct routes and the second $N^\circ$ queues corresponding to the number of clients on the indirect routes. The parameters $\gamma_n^{({\rm u})}$ and $\gamma_n^{({\rm d})}$ correspond to the up-time and down-time of link $n$. Clients for origin-destination pair $n$ arrive according to a Poisson process with rate $\lambda_n$, and stay in the system for an exponential time with parameter $\mu_n.$ We again stress that various generalizations are possible, such as phase-type up- and down-times and Markov modulated arrival processes; these extensions are conceptually very similar to the setup we describe here, but notationally burdensome.

Each of the $N^\circ$ links can be up or down, so that the background process has $I=2^{N^\circ}$ states. Suppose the background process is in state $i$. Again, $1_n(i):=1\{n \in S(i)\}$.
For $n=1,\ldots,N^\circ$,
\[\lambda_n^{(i)} = \lambda_n\,1_n(i)
,\:\:\: \lambda^{(i)}_{n+N^\circ} = \lambda_n \,(1-1_n(i))\,1_{n_1(n)}(i)\, 1_{n_2(n)}(i).\]
All $\mu_{nn'}=0$ for $n,n'\in\{1,\ldots,N\}$, and $\mu_{n0}= \mu_{n+N^\circ,0}=\mu_n.$

We now consider the  transitions corresponding to links going down (and coming up again). We distinguish two cases; for all $i,j\in{\mathscr I}$ we have that $K_{ij}$ equals 0 or 1. 
\begin{itemize}
\item[$\circ$] Suppose that for a $j\not= i$ and some $n\in\{1,\ldots,N^\circ\}$   we have $S(i)=S(j)\cup\{n\}$; then the background process jumps from $i$ to $j$ after an exponentially distributed time with rate $\alpha_{ij}=\gamma_n^{({\rm u})}$. Note that this transition corresponds to link $n$ failing, and thus clients using this route as a direct route move to the indirect route (if available) and
clients using this link as part of their indirect route are lost.
The  queue content vector is premultiplied by a $(N\times N)$-dimensional matrix $A_{ij}$ consisting of (i)~a $0$ on position $(n,n)$, (ii)~a $1$ on position $(n+N^\circ,n)$ but only if ${\mathscr N}(n)\subseteq S(i)$ (where it is noted that if
${\mathscr N}(n)\not\subseteq S(i)$, then files are lost), (iii)~a $0$ on position $(n'+N^\circ,n'+N^\circ)$ if $\{n\}\subseteq {\mathscr N}(n')$ (corresponding to files that are lost),  (iv)
all other diagonal entries  being 1, and  (v) all other entries being~0.
\item[$\circ$] 
Suppose on the other hand that for $i\not=j$ and some $n\in\{1,\ldots,N^\circ\}$ we have $S(j)=S(i)\cup\{n\}$; then the background process jumps from $i$ to $j$ with rate $\alpha_{ij}=\gamma_n^{({\rm d})}$. This transition corresponds to link $n$ having been repaired. The  queue content vector is premultiplied by a $(N\times N)$-dimensional matrix $A_{ij}$ consisting of (i) a $0$ on position $(n+N^\circ,n+N^\circ)$, (ii) a~$1$ on position $(n,n+N^\circ)$,
(iii)
all other diagonal entries  being 1, and  (iv) all other entries being~0.
\end{itemize}

Again, our model can be used to study design questions.
As indicated above, clients are lost if both the direct route and the indirect route are unavailable. Compared to using only direct routes, the option of indirect routes evidently reduces the number of lost clients, but this comes at the price of the servers being more intensively used.  Let $Z_\ell(t)$ denote, as before, the number of clients lost in $[0,t]$; see Equation (\ref{LOST}). In addition, let $Z_s(t)$ be the amount of link resources used in $[0,t]$:
\[{\mathbb E}\,Z_s(t) = \int_0^t \sum_{i=1}^I \sum_{n=1}^{N^\circ} {\mathbb E}\big(M_n(s) I_i(s)\big)\,{\rm d}s
+2  \int_0^t \sum_{i=1}^I \sum_{n=N^\circ+1}^{N} {\mathbb E}\big(M_n(s) I_i(s)\big)\,{\rm d}s
;\]
again Lemma \ref{lemma2} can be used to numerically evaluate this quantity.

With a time horizon $T$, let $\varrho_\ell$ be the cost of loss and $\varrho_s$ the cost of storage, so that the total cost is
\begin{equation}
\label{obj}
\varrho_\ell \cdot {\mathbb E}\,Z_\ell(T) + \varrho_s \cdot {\mathbb E}\,Z_s(T).\end{equation}
Its value can be compared to the value of the same objective function in the system without rerouting. Typically, for small $\varrho:=\varrho_\ell/\varrho_s$ the system without rerouting is to be preferred, whereas for large $\varrho$ rerouting pays off. To optimally design the system, it would be useful to have knowledge of the critical value $\varrho\s$ (for which both mechanisms have the same cost, that is). 

\subsection{Applications to storage networks}\label{sto}
In storage networks information is typically stored at multiple locations (e.g.\ on multiple data storage servers), so as 
to mitigate the danger of files being lost. A relevant design issue concerns developing a policy 
that controls the fraction of files lost without unnecessarily replicating them. 
For a general account of various aspects of storage networks, see e.g.\ 
\cite{PES}.

Consider a system with $K$ storage locations, each of which  can be either  `up' or `down'.
Let the up-time of location $k\in\{1,\ldots,K\}$ be exponentially distributed with parameter $\gamma_k^{({\rm u})}$, and let the corresponding down-time exponentially distributed with parameter $\gamma_k^{({\rm d})}$. We  thus have constructed an environmental process of dimension $I=2^K$, where each state corresponds to the set of locations that are up (and consequently also the set of locations that are down). We let, for any $i\in\{1,\ldots,I\}$, the set $U(i)$ denote the locations that are up when the environmental process is in state $i$. We order the $I$ states such that the state $1$ corresponds to all locations up, the states $2$ up to $K+1$ to all situations with one location down, etc., so that state $2^K$ corresponds to all locations being down.

Files can be stored on any subset of the locations; there are $N=2^K-1$ of these. We let $S(n)$ denote the locations involved in the $n$-th subset, for $n\in\{1,\ldots,N\}$. These are ordered in the same way as above: 
queue $1$ corresponds to files stored at all locations, the queues $2$ up to $K+1$ to files stored at all-but-one locations, etc., so that queues  $2^K-K-2$ up to $2^K-1$ correspond to files stored on just a single location (which are lost if this location fails).

We now argue that this model is covered by the general multiplicative-transition framework that we introduced in the previous section. 
Consider the situation that the environmental process is in state $i$. 
Let $\lambda_n$ be the (constant) arrival rate that is intended to be stored at the set of locations $S(n)$. However, if $i$ is such that this is not possible (because some of the locations are down), it is only stored at the subset of $S(n)$ that is up.
This means that, with $V(i,n):=\{n': S(n')\cap U(i) = S(n)\}$, external arrivals to subset $n$ occur at rate
\[\lambda_n^{(i)} = \sum_{n'\in V(i,n)} \lambda_{n'}.\]
During operations, files may be copied to additional locations, may be deleted from locations or may be deleted completely. Therefore, files hop
from queue $n$ to $n'$ with rate $\mu_{nn'}^{(i)}$ (with $n'=0$ corresponding to files leaving the network).  

We now consider the multiplicative transitions. Two cases are to be distinguished. 
\begin{itemize}
\item[$\circ$]
Suppose that for some $j\in\{1,\ldots,I\}$, that is assumed to be different from the current environmental state $i$,  it holds that $U(i)=U(j)\cup\{k\}$; then the background process jumps from $i$ to $j$ after an exponentially distributed time with rate $\gamma_k^{({\rm u})}$ (note that this transition corresponds to the event that location $k$ finishes its up-time, i.e., goes down). Simultaneously the $N$-dimensional queue content vector is premultiplied by a matrix $A_{ij}$ that is defined as follows. It has a zero on the diagonal positions that correspond to subsets that contain location $k$ (i.e., $n$ such that $\{k\}\subseteq S(n)$).  In the same column, it has a one on the position $n'$ such that $S(n')= S(n)\setminus\{k\}$ (if any).
\item[$\circ$]
Suppose that for $i\not=j$ we have $U(j)=U(i)\cup\{k\}$; then the background process jumps from $i$ to $j$ with rate $\gamma_k^{({\rm d})}$ (without any change in the network population vector; this transition corresponds to the event that location $k$ finishes its down-time, i.e., becomes functioning again). 
\end{itemize}
Recalling that the entries  $2^K-K-2$ up to $2^K-1$ of ${\bs M}(\cdot)$ correspond to files stored at just a single location, 
we can evaluate the mean number of lost files in $[0,t]$   as
\[{\mathbb E}\,Z_\ell(t) = \int_0^t\sum_{i=1}^I\sum_{n=2^K-K-2}^{2^K-1}
{\mathbb E}\big(M_n(s)I_i(s)\big) \gamma^{\rm (u)}_{n-2^K+K+1} {\rm d}s,\]
which can be numerically evaluated using Lemma \ref{lemma2}.

\vb

Consider for example the case of $K=2$ locations, so that $I=4$ and $N=3$. In self-evident notation we code the $4$ states of the background process as 
\[\{1,2,3,4\}\equiv\{\{1,2\},\{1\},\{2\},\emptyset\}\]
(with the left-hand side in the previous display being in terms of the elements $i\in{\mathscr I}$, and the  right-hand side in terms of the corresponding $U(i)$).
Then $K_{ij}=1$ for all $i,j\in{\mathscr I}$, and
\[A_{12}=A_{34}=\left(\begin{array}{ccc}
0&0&0\\ 1&1&0\\ 0&0&0\end{array}\right),\:\:
A_{13}=A_{24}=\left(\begin{array}{ccc}
0&0&0\\ 0&0&0\\ 1&0&1\end{array}\right),
\]
whereas the other $A$-matrices equal $I_3$ (note that $A_{12}$ and $A_{34}$ correspond to location 2 going down, and $A_{13}$ and $A_{24}$ to location 1 going down). In addition, 
\[\alpha_{12}=\alpha_{34}=\gamma_2^{({\rm u})},\:\:\alpha_{13}=\alpha_{24}=\gamma_1^{({\rm u})},\:\:
\alpha_{21}=\alpha_{43}=\gamma_2^{({\rm d})},\:\:\alpha_{31}=\alpha_{42}=\gamma_1^{({\rm d})}.\]

\section{Numerical experiments}
To illustrate the potential of our results,  in this section we provide two examples: one on a retrial queue, and another one on storage networks.

\subsection{Retrial queue}
{In this first example, we consider a single retrial system, i.e., a two-queue network consisting of a service station and a retrial location. The service station alternates between being `up' and `down', with the corresponding durations being exponentially distributed with parameters $\gamma^{\rm (u)}$ and $\gamma^{\rm (d)}$, respectively. Clients arrive with rate $\lambda$ and their service times are exponentially distributed with mean $\mu^{-1}$. The rate at which customers in the retrial queue attempt to reenter service is $\kappa$, where the corresponding renege rate is $\nu.$

We now cast this system in the terminology of our overarching model. 
The background process can be in two states (so that $I=2$); we let state $1$ correspond to the station being functioning, and state $2$ to the station being down. The dimension of ${\bs M}(\cdot)$ is $N=2$; the first component corresponds to the queue at the service station, whereas the second component corresponds to the retrial pool. The matrices $A_{12}$ and $A_{21}$ are given by
\[A_{12} = \left(\begin{array}{cc}0&0\\
1&1\end{array}\right),\:\:\:\:A_{21} = \left(\begin{array}{cc}1&0\\
0&1\end{array}\right).\] The arrival rates are $\lambda_n^{(i)}=\lambda$ for $(i,n)$ equalling $(1,1)$ or $(2,2)$, and otherwise 0. 
In addition, $\mu_{21}^{(1)} = \kappa$, $\mu_{20}^{(1)}=\mu_{20}^{(2)} = \nu$, $\mu_{10}^{(1)}=\mu$, whereas the other departure rates are $0$. Also, $\alpha_{12}= \gamma^{\rm (u)}$ and $\alpha_{21}= \gamma^{\rm (d)}$.

Let $\bar M_{ni}(t)$ be the mean number in queue $n$ when the background process is in state $i$ at time $t$; observe that we constructed our model such that $\bar M_{12}(t)=0$ for all $t\geqslant 0$. The time-dependent means follow from solving a system of linear differential equations:
\begin{align*}
\bar M'_{11}(t) &= \lambda \pi_1(t) - (\mu+\gamma^{\rm (u)}) \bar M_{11}(t) + \kappa \bar M_{21}(t),\\
\bar M'_{21}(t) &= \gamma^{\rm (d)} \bar M_{22}(t)-(\kappa + \nu+\gamma^{\rm (u)}) \bar M_{21}(t),\\
\bar M'_{22}(t) &= \lambda\pi_2(t)+ \gamma^{\rm (u)}\bar M_{11}(t)+ \gamma^{\rm (u)}\bar M_{21}(t)-(\nu+
\gamma^{\rm (d)})
 \bar M_{22}(t).
\end{align*}
We now present the stationary means $\bar M_{11}$, $\bar M_{21}$, and $\bar M_{22}$.
Let $\Gamma:= \gamma^{\rm (u)}+\gamma^{\rm (d)}$, $\pi_1 = \gamma^{\rm (d)}/\Gamma = 1-\pi_2.$ Sending $t\to\infty$, and letting the derivatives in the above differential equations be equal to 0, we obtain
\[\bar M_{21} = \frac{\lambda\, \gamma^{\rm (u)}}{\Gamma\eta}\left(\frac{\mu+\gamma^{\rm (u)}+\gamma^{\rm (d)}}{\mu+\gamma^{\rm (u)}}\right),\:\:\:
\eta:=(\kappa+\nu+\gamma^{\rm (u)})\frac{\nu+\gamma^{\rm (d)}}{\gamma^{\rm (d)}}-\kappa \frac{\gamma^{\rm (u)}}{\mu+\gamma^{\rm (u)}}-\gamma^{\rm (u)},\]
and
\[\bar M_{11} = \frac{1}{\mu+\gamma^{\rm (u)}}\left(\kappa \bar M_{21}+\lambda
\frac{ \gamma^{({\rm d})}}{\Gamma}\right),\:\:\:
\bar M_{22} =\frac{\kappa+\nu+\gamma^{\rm (u)}}{\gamma^{\rm (d)}} \bar M_{21} 
.\]
We now consider the model's loss ratio $\ell$, defined as the long-run fraction of clients leaving the network without being served (i.e., due to reneging). With $\bar M_{21}$ and $\bar M_{22}$ as computed above, 
\[\ell = \frac{\nu}{\lambda}\big(\bar M_{21}+\bar M_{22}\big).\]

{\it Experiment 1.} To control the loss ratio, the service provider may opt for speeding up the repair times. The above formulas allow us to determine the smallest $\gamma^{\rm (d)}$ such that the loss ratio $\ell$ is below some maximally allowed value $\ell\s$; observe that $\ell$ is decreasing in $\gamma^{\rm (d)}$. It requires elementary calculus to verify that
\[\lim_{\gamma^{\rm (d)}\to\infty} \bar M_{21} = \frac{\lambda\gamma^{({\rm u})}}
{\kappa\mu+\nu\mu +\nu \gamma^{\rm (u)}},\:\:\:\lim_{\gamma^{\rm (d)}\to\infty} \bar M_{22} =0,\]
so that
\[\lim_{\gamma^{\rm (d)}\to\infty} \ell=\frac{\nu\gamma^{({\rm u})}}
{\kappa\mu+\nu\mu + \nu\gamma^{\rm (u)}};\]
this expression increases in $\gamma^{\rm (u)}$ and $\nu$ and decreases in $\mu$ and $\kappa$, as expected.

Observe that it in general cannot be guaranteed that there is a $\gamma^{\rm (d)}$ such that $\ell\leqslant \ell\s$: the parameters can be such that $\ell>\ell\s$ for all $\gamma^{\rm (d)}$.
This is because even very short down-times lead to the event of clients simultaneously moving to the retrial queue, where the effect of clients reneging starts to kick in.


In the numerical experiment  we chose $\lambda=100$, $\kappa=2$, $\nu=2$, $\mu=1$ and $\gamma^{({\rm u})}=0.1$. First suppose that the loss ratio should remain below $10\%$. 
One needs to   take $\gamma^{({\rm d})}$ larger than $2.1496$, as illustrated by Fig.\ \ref{F1a} (left panel).
Suppose, on the contrary, that the target is $1\%$, then this cannot be achieved by increasing $\gamma^{({\rm d})}$; based on the above results, we conclude that even by making the repairs very fast, the loss ratio will (for these values of $\lambda$, $\kappa$, $\nu$, $\mu$ and $\gamma^{({\rm u})}$) never get below $0.2/4.2 \approx 4.76\%$ (corresponding to the horizontal dashed line in the graph). 

\begin{figure}
\includegraphics[scale=1.090]{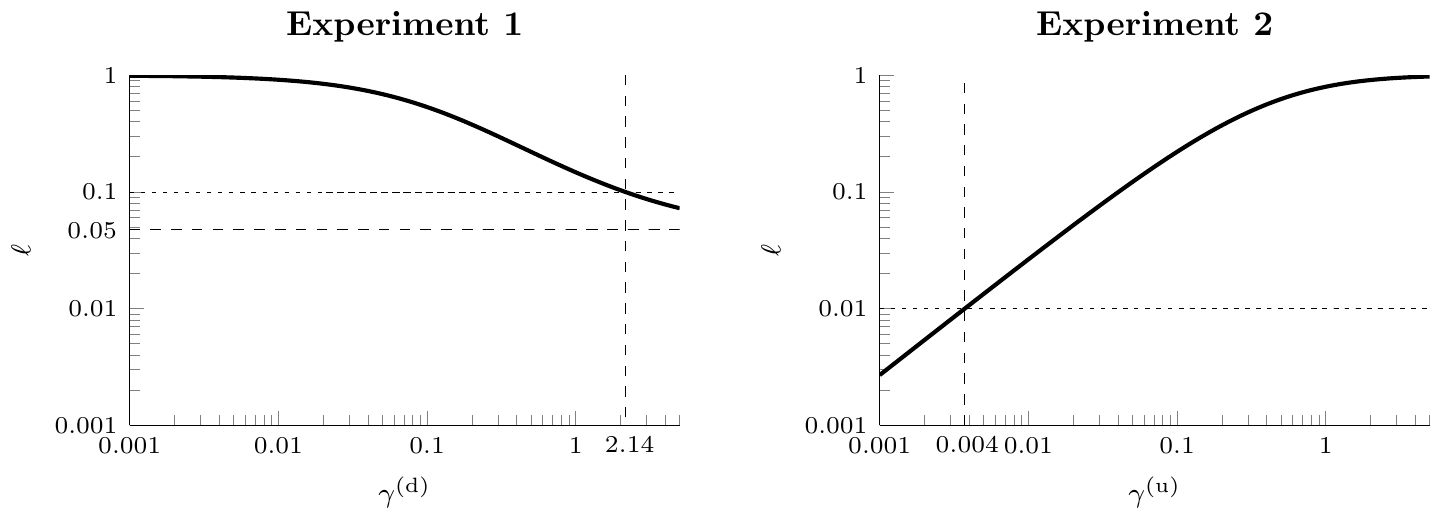}
\caption{\label{F1a} Retrial queue: loss ratio $\ell$, Experiments 1 and 2.
}
\end{figure}

\vb

{\it Experiment 2.} An alternative way to control $\ell$ is by making the up-times longer, i.e., by decreasing $\gamma^{\rm (u)}$. It is readily verified that
\[\lim_{\gamma^{\rm (u)}\downarrow 0} \bar M_{21} = \lim_{\gamma^{\rm (u)}\downarrow 0} \bar M_{22} =0,\]
so that the loss rate $\ell$ will be below any critical value $\ell\s$ for $\gamma^{\rm (u)}$ small enough. 

In our numerical experiment we again chose   $\lambda=100$, $\kappa=2$, $\nu=2$, $\mu=1$, but now we fix $\gamma^{({\rm d})}=0.5$. We wonder whether in this scenario a loss ratio   below $1\%$ can be achieved by tuning $\gamma^{({\rm u})}$. Fig.\ \ref{F1} (right panel) shows that this is indeed the case: as it turns out, $\gamma^{({\rm u})}$ should be below  $0.0037.$ 

\vb

In practice, one may want to find the most cost effective pair $(\gamma^{\rm (u)},\gamma^{\rm (d)})$ such that the performance requirement is met. 
With $\varrho^{\rm (u)}$  the cost of making the mean up-times one  unit longer, and $\varrho^{\rm (d)}$
the cost of making the hazard rate corresponding to the down-times one unit larger, a relevant optimization problem could be of the type
\[\min_{\gamma^{\rm (u)},\gamma^{\rm (d)}} \frac{\varrho^{\rm (u)}}{\gamma^{\rm (u)}}+\varrho^{\rm (d)} \gamma^{\rm (d)},\:\:\:\mbox{subject to}\:\:\ell\leqslant \ell\s.\]

\subsection{A storage system}{   
In this example we show how to analyze a specific storage system; it has some elements   
in common with the class of models that was introduced in Section \ref{sto}, but there 
are a few notable differences. Files arrive according to a Poisson process with rate $\lambda$. 
With probability $p$ the file will be offered standard service, and with probability $1-p$ 
premium service. A {\it basic file} is randomly allocated to one of the two locations (say A and B), 
where it will stay until that location goes `down'. In our example copies of files are never 
deleted (except through a failure of the storage location).
A {\it premium file} is randomly allocated to one of the locations, but is then copied at rate 
$\mu$ to the other location.
When a location goes `down' in the premium case, upon repair files are again copied back 
(also at rate $\mu$, that is). The locations' up- and down-times are independent and 
statistically identical; up-times (down-times, respectively) are exponentially distributed 
with rate $\gamma^{\rm (u)}$ ($\gamma^{\rm (d)}$). 
In this system there are five queues to be kept track of: premium files on location A, 
premium files on location B, premium files on locations A and B, basic files on location  A, 
and basic files on location B.

\vb

{\it Experiment 1.} The parameters we picked are: $\lambda = 10\,000$ (i.e., on average 10\,000 files arrive per day), $\mu = 24$ (i.e., it takes on average an hour before a stored file is copied to a second location), $\gamma^{\rm (u)}=0.01$ (i.e., each of the storage locations are functional on average for consecutive periods of 100 days), and $\gamma^{\rm (d)}=2$ (i.e., it takes 12 hours to repair a storage location). We let the system start empty at time 0, with both locations being `up' (but other initial conditions are handled in the precise same way). 

The first graphs show, for $T=1$ (i.e., one day), the expected number of lost files ${\mathbb E}\,Z_\ell(T)$, and the expected integral of the number of stored files,    ${\mathbb E}\,Z_s(T)$, as functions of the fraction of premium files $p$. In the previous section we pointed out how these metrics can be evaluated, but the computation of ${\mathbb E}\,Z_\ell(T)$ can be done more efficiently, relying on the following idea; the performance measure ${\mathbb E}\,Z_a(T)$ can be dealt with analogously.

The idea is to append one coordinate to the state space; the resulting extra component  $M_{N+1}(t)$ records the number of files lost in $[0,t]$ (which can be seen as a queue with zero departure rate). The  transform of the vector $({\bs M}(t),M_{N+1}(t))\in{\mathbb N}_0^{N+1}$ (jointly with the state of the environmental process) can be characterized in the precise same way as that of just ${\bs M}(t)$, i.e., by setting up a system of partial differential equations. This provides us with an expression for ${\mathbb E}\,Z_\ell(T)$
of the form~(\ref{MEM}). Observe that it entails that we can evaluate the quantity ${\mathbb E}\,Z_\ell(T)$, which can be evaluated using Lemma \ref{lemma1}; in this way we  
avoid  evaluating integrals of the type of (\ref{LOST}).

\begin{figure}
\includegraphics[scale=1.090]{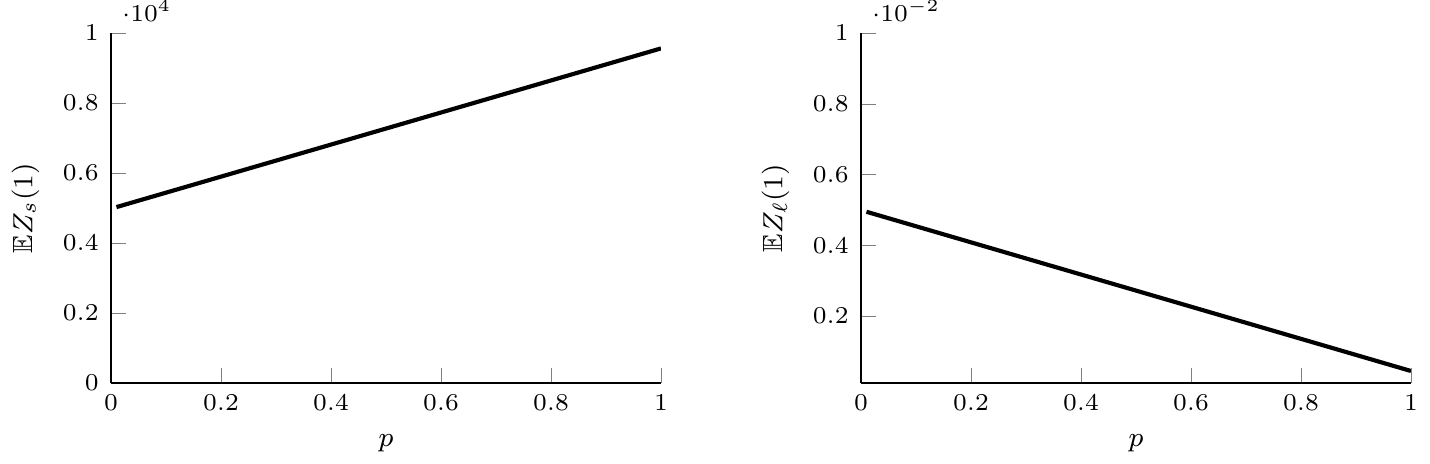}
\caption{\label{F1} Storage system: ${\mathbb E}\,Z_\ell(T)$ and ${\mathbb E}\,Z_s(T)$ as functions of $p$, Experiment 1.
}
\end{figure}

The graphs in Fig.\ \ref{F1} show, for $T=1$, that 
${\mathbb E}\,Z_s(T)$ increases in $p$ (left panel), whereas 
${\mathbb E}\,Z_\ell(T)$ decreases in $p$ (right panel),  as expected.

\vb

{\it Experiment 2.}    We now consider a cost function that is a linear combination of ${\mathbb E}Z_\ell(t)$ and ${\mathbb E}Z_s(t)$, i.e., (7). In this case  the optimal design amounts to minimizing the objective function (7) with respect to the fraction $p \in [0,\,1]$. Let $\varrho_\ell$ and $\varrho_s$ again respectively correspond to the cost of a lost file and the cost of a unit of storage per unit time.  Clearly,  $p^\star = 0$ for $\varrho \downarrow 0$ (as losing files is not penalized), whereas $p^\star = 1$ for $\varrho \uparrow \infty$ (as storing files is  not penalized). Bearing in mind the  shapes  of ${\mathbb E}Z_\ell(t)$ and ${\mathbb E}Z_s(t)$, as depicted in Fig.\ \ref{F1}, the optimization of a linear combination of ${\mathbb E}Z_\ell(t)$ and ${\mathbb E}Z_s(t)$ leads to 
$p^\star$ equalling either  0 or 1. 
The left panel of Fig.~\ref{F2} shows the region in which the optimal $p^\star$ is 0 or 1, for combinations of 
$\gamma^{({\rm d})}$ and $\varrho := \varrho_\ell/\varrho_s$ with $\varrho_s$ fixed equal to one and $\gamma^{({\rm u})}$ equal to one (and all other parameters as in Experiment 1).
In the right panel of Fig.~\ref{F2} we show a similar picture, but now with $\gamma^{({\rm u})}$ on the horizontal axis.

\begin{figure}
\includegraphics[scale=1.090]{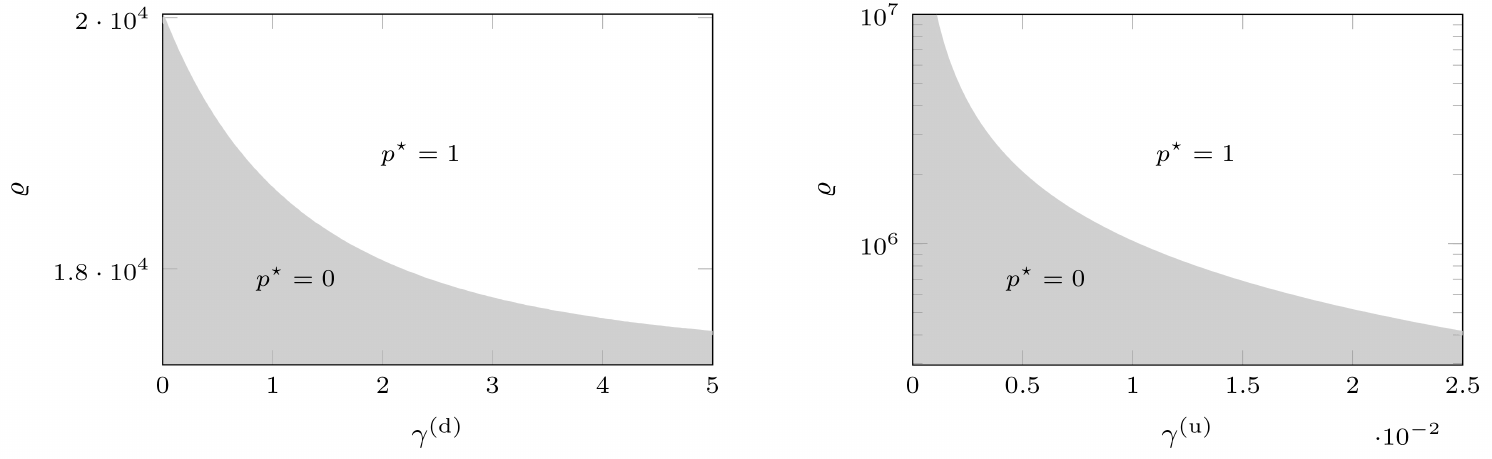}
\caption{\label{F2} Storage system: areas in which $p^\star$ equals 0 and 1, for different values of the rates $\gamma^{({\rm d})}$ and 
$\gamma^{({\rm u})}$ on the horizontal axis and `cost ratio' $\varrho$ on the vertical axis, Experiment 2.
}
\end{figure}

\begin{figure}[h]
	\centering
  \includegraphics{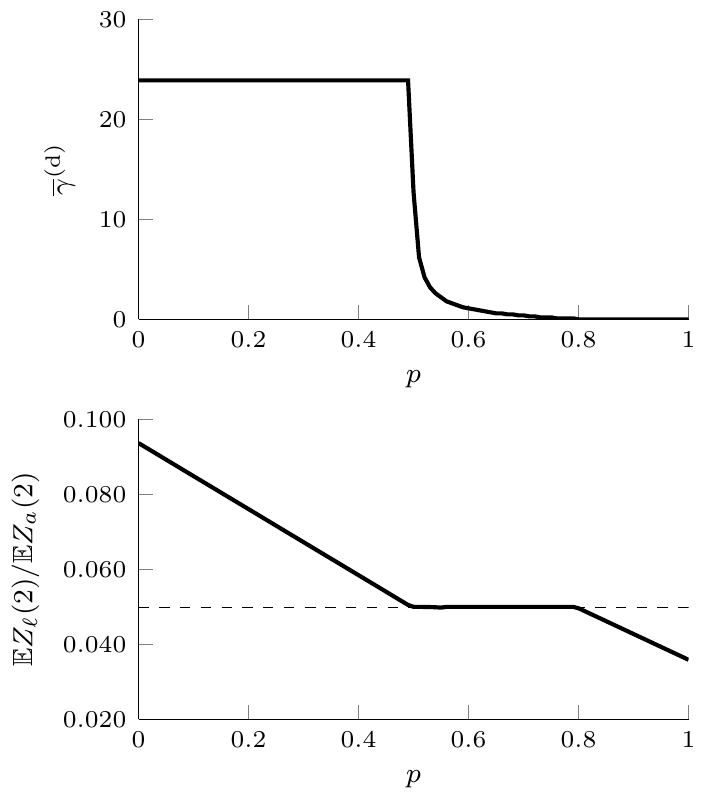}
	\caption{Storage system: $\gamma^{{\text (d)}}$ and the corresponding proportion of files lost for different values of the fraction $p$, Experiment 3. }
	\label{fig:VaryP}
\end{figure}\vb

{\em Experiment 3.} We now vary the value of the repair rate $\gamma^{({\rm d})}$ with the goal of achieving a predetermined performance target. For any value of $p$ we compute the minimally required repair rate (defined as $\bar\gamma^{({\rm d})}$) from $\gamma^{({\rm d})} \in [0, 24]$, in an attempt to ensure that the loss fraction ${\mathbb E}Z_\ell(T)/{\mathbb E}Z_a(T)$ is below 0.05 (where we pick $T=2$). Observe that the constraint $\bar\gamma^{({\rm d})}\leqslant 24$ amounts to imposing the 
requirement that repairs must be expected take at least 1 hour to perform.

Inspection of Fig.\ \ref{fig:VaryP} immediately reveals that for $p$ smaller than 0.5 we are unable to achieve our desired loss fraction using only the available changes in $\gamma^{({\rm d})}$. Indeed it is conceivable that for small $p$ there does not exist a repair rate such that the loss fraction goes below $0.05$, a phenomenon similar to that which we earlier saw in Experiment~1 for the retrial queue. As $p$ is increased within $[0,0.5]$ we see an approximately linear decrease in the loss fraction resulting from the increased proportion of files being placed in the premium category (where they are unlikely to become lost).  For $p \in [0.5,0.8]$, we observe  that we are able to achieve our desired loss fraction; moreover the storage location can be repaired increasingly slowly if more files are multiply stored (i.e., when $p$ increases). This effect initially results in a very rapid decrease in the repair rate but has less of an impact as $p$ is increased closer to 0.8, at which point the repair rate can no longer be traded off against increased duplication. Notice that the mechanism by which $\gamma^{({\rm d})}$ decreases basic file losses is by reducing the portion of $[0,T]$ during which {\em both} storage locations are inoperable; this variable has no effect on basic files which are accepted into the system only to be lost due to a  failure later. Hence, focusing on basic files, it can be seen that eventually the effect of $\gamma^{({\rm d})}$ on the portion of $[0,T]$ for which both storage locations  are inoperable becomes negligible compared to the reductions in losses from increasing $p$. The result of this is that for $p>0.8$ the loss fraction continues to decrease approximately linearly as more files are placed in the very safe premium category, as we saw for $p<0.5$.

\section{Discussion and concluding remarks}
In this paper we studied a network of Markov-modulated infinite-server queues with the distinguishing feature that it also incorporates events by which the network population vector makes {\it multiplicative transitions} (at which it changes from 
${\bs m}$ to $A{\bs m}$, for some matrix $A$). As we argued, the resulting framework covers various relevant models as special cases; for example, it enables the modelling of retrial queues, networks with rerouting, and storage systems. 

Our results for the system's transient behavior are in terms of (i)~a system of partial differential equations describing the moment generating function of the network population vector, and (ii)~a procedure to compute moments. In these expressions time $t$ can be sent to $\infty$ so as to obtain the corresponding stationary behavior, under the proviso that the stability condition applies. 

\vb

{\it Future research.} The model we have developed triggers various intriguing research questions. In the first place, one may wonder whether under a specific scaling of the parameters one could find a weak limit for its transient or stationary behavior. Such a procedure has been developed in \cite{BKM, KOEN} for Markov-modulated infinite-server queues {\it without} multiplicative transitions. For that model the limiting process (after scaling the arrival rates and the environmental process) is a multivariate Ornstein-Uhlenbeck process. In this diffusion limit all limiting marginal distributions (and the model's stationary distribution, too) are asymptotically Normal. 
For our model however, {\it with} multiplicative transitions that is,  it is anticipated that there is no limiting process of diffusion type, due to the possibly large jumps caused by the multiplicative transitions; cf.\ \cite{FR}. In particular, the marginal distributions are expected to be asymmetric. 

Scaling the external arrival rates by a common factor, say $K$, it is seen from (\ref{statm}) that the stationary mean also grows proportionally to $K$. Calling the stationary distribution under this scaling ${\bs M}^{(K)}$, one may want to asymptotically characterize large-deviation probabilities of the type
\[p_K:={\mathbb P}\left(\frac{{\bs M}^{(K)}}{K}\in S\right),\]
for $K$ large and  a set $S$ that does not contain ${\mathbb E}\,{\bs M}^{(K)}/K=-({\mathscr M}+{\mathscr A})^{-1} {\mathscr L}  {\bs \pi}$. 
It is not clear how such asymptotics  can be found; observe that due to the multiplicative transitions the model does not fit in the Freidlin-Wentzell framework \cite{ShW}, so that standard large-deviation techniques are likely to fail.

Other challenges lie in the application of our techniques to develop design principles for various sorts of operational networks. For instance for storage networks, one may want to develop an optimal replication policy, striking a proper balance between controlling the risk of files being lost and excessively using storage space. 

{\small
\subsection*{Acknowledgments}

The authors wish to thank Peter Taylor (The University of Melbourne) and Ross McVinnish (The University of Queensland) for useful remarks.}

\bibliographystyle{plain}

\end{document}